%% file: nonlocal-obst-survey-180727.tex
\documentclass[11pt,reqno]{amsart}

\input{latexformat}
\input{latexnewmacros}

\numberwithin{equation}{section}
\theoremstyle{default}

\newtheorem{thrm}{Theorem}[section]
\newtheorem{dfn}[thrm]{Definition}
\newtheorem{lemma}[thrm]{Lemma}
\newtheorem{prp}[thrm]{Proposition}

\usepackage{hyperref}
\usepackage{mathrsfs}
\usepackage{slashed}
\usepackage[usenames]{color}
\usepackage{url}
\usepackage{verbatim}
\usepackage{amssymb,amsmath, graphicx}
\usepackage{amsfonts}
\usepackage{amsthm}
\usepackage{enumerate, latexsym, esint}
\usepackage{amscd, amsmath}
\usepackage{tikz}

\allowdisplaybreaks

\hypersetup{pdftitle={Nonlocal operators}}

\hypersetup{pdfauthor={Donatella Danielli, Arshak Petrosyan, and Camelia A. Pop}}

\begin{document}

\title[Obstacle problems for nonlocal operators]{Obstacle problems for nonlocal operators: A brief overview}

\author[D. Danielli]{Donatella Danielli}
\address[DD]{Department of Mathematics, Purdue University, West Lafayette, IN 47907}
\email{danielli@math.purdue.edu}

\author[A. Petrosyan]{Arshak Petrosyan}
\address[AP]{Department of Mathematics, Purdue University, West Lafayette, IN 47907}
\email{arshak@math.purdue.edu}

\author[C. A. Pop]{Camelia A. Pop}
\address[CP]{School of Mathematics, University of Minnesota, Minneapolis, MN 55455}
\email{capop@umn.edu}

\begin{abstract}
In this note, we give a brief overview of obstacle problems for nonlocal operators, focusing on the applications to financial mathematics. The class of nonlocal operators that we consider can be viewed as infinitesimal generators of non-Gaussian asset price models, such as Variance Gamma Processes and Regular L\'evy Processes of Exponential type. In this context, we analyze the existence, uniqueness and regularity of viscosity solutions to obstacle problems which correspond to prices of perpetual and finite expiry American options.
\end{abstract}

%\date{\today{ }\hhmm}

\maketitle

%\tableofcontents

\section{Introduction}
\label{sec:Introduction}

The purpose of this note is to give a brief overview of obstacle problems for nonlocal operators, focusing on the applications to financial mathematics. Natural classes of nonlocal operators are infinitesimal generators of \emph{L\'evy processes.}
We recall that a L\'evy process $\{X(t)\}_{t\geq 0}$
defined on a filtered probability space $(\Omega,\mathcal{F}, \{\mathcal{F}_t\}_{t\geq 0}, \PP)$
is a random process that is stochastically continuous and has stationary and independent increments.
More precisely, $\{X(t)\}_{t\geq 0}$  is a L\'evy process if:
\begin{itemize}
\item[1.] $X(0)=0$ with probability 1;
\item[2.] For all $0\leq t_1<t_2<\dots<t_n$,  $X(t_1),\ X(t_2)-X(t_1), \dots,\ X(t_n)-X(t_{n-1})$ are independent;
\item[3.] For all $0\leq s<t<\infty$, the probability distribution of $X(t)-X(s)$ is the same as the one of $X(t-s)$;
\item[4.] For all $\eps>0$, we have that
$$
\lim_{t\downarrow 0} \PP\left(|X(t)| > \eps\right) = 0.
$$
\end{itemize}
We begin the introduction with \S\ref{sec:Rep_Levy_proc} where we gives representations of L\'evy processes using the L\'evy-Khintchine formula and the L\'evy-It\^o decomposition. We continue in \S\ref{sec:Connection_to_pde} to describe the connection to nonlocal (integro-differential) operators and we present in \S\ref{sec:SDE} more general stochastic equations, which give rise to a wider class of nonlocal operators. Finally, in \S\ref{sec:Obs_prob} we give a brief introduction to obstacle problems and we summarize in \S\ref{sec:Literature_overview} previous results obtained in the literature.

\subsection{Representations of L\'evy processes}
\label{sec:Rep_Levy_proc}
Our starting point is the \emph{L\'evy-Khintchine formula} \cite[Corollary~2.4.20]{Applebaum}, which shows that, for all $t\geq 0$ and $\xi\in\mathbb{R}^n$, we have
\begin{equation}
\label{eq:LK}
\mathbb{E}\left[e^{i\xi\cdot X(t)}\right]=e^{t\psi(\xi)},
\end{equation}
where the characteristic exponent $\psi(\xi)$ is given by
\begin{equation}
\label{eq:exponent}
\psi(\xi)=-\frac{1}{2}\xi\cdot A\xi + ib\cdot \xi+\int_{\mathbb{R}^n\setminus\{0\}}\left(e^{i\xi\cdot y}-1-i\xi\cdot y\chi_{|y|<1}\right)\ \nu(dy).
\end{equation}
Here $A$ is a $n\times n$-dimensional, symmetric, positive-semidefinite matrix, $b\in\mathbb{R}^n$ and $\nu$ is a L\'evy measure on $\mathbb{R}^n\setminus\{0\}$, i.e. it satisfies
$$
\int_{\mathbb{R}^n\setminus\{0\}} \min\{1,|y|^2\}\ \nu(dy)<\infty.
$$
When $A\equiv 0$ and $\nu\equiv 0$, that is $\mathbb{E}\left[e^{i\xi\cdot X(t)}\right]=e^{itb\cdot \xi}$,  the process $X(t)=tb$ is deterministic motion on a straight line, with velocity of motion, or \emph{drift}, $b$.
If instead $A\equiv 0$, but $\nu\not\equiv 0$ has finite variation, that is it satisfies
\begin{equation}
\label{eq:finite_variation}
\int_{\mathbb{R}^n\setminus\{0\}} \min\{1, |y|\}\ \nu(dy) < \infty,
\end{equation}
%removed extra curly bracket in domain of integration
then we can rewrite the characteristic exponent \eqref{eq:exponent} as
$$
\psi(\xi)=ib'\cdot \xi+\int_{\mathbb{R}^n\setminus\{0\}}\left(e^{i\xi\cdot y}-1\right)\ \nu(dy).
$$
%Maybe going through these examples can be shortened
The simplest possible case is when $\nu=\lambda \delta_h,$ where $\lambda>0$ and $\delta_h$ is the Dirac mass concentrated at $h\in\mathbb{R}^n\setminus\{0\}$. If we let $X(t)=b't+N(t)$, then the process $\{N(t)\}_{t\geq 0}$ is such that
$$
\mathbb{E}\left[e^{i\xi\cdot N(t)}\right] = \exp{\left[\lambda t\left(e^{i\xi\cdot h}-1\right)\right]},
$$
and therefore $\{N(t)\}_{t\geq 0}$ is a Poisson process of intensity $\lambda$ taking values in $\{mh,\ m\in\mathbb{N}\}$. The physical interpretation is that $\{X(t)\}_{t\geq 0}$ follows the path of a straight line with drift $b'$ and has jump discontinuities of size $|h|$. The time between two consecutive jumps are independent random variables exponentially distributed with parameter $\lambda$.

The next step is to take $\nu=\sum_{j=1}^m \lambda_j\delta_{h_j}$, with $m\in\mathbb{N}$, $\lambda_j>0,\ h_j\in\mathbb{R}^n\setminus\{0\},\ 1\leq j\leq m.$ In this instance, we can write $\{X(t)\}_{t\geq 0}$ as
$$
X(t)=b't+\sum_{j=1}^m N_j(t),
$$
where the $\{N_j(t)\}_{t\geq 0},\ 1\leq j\leq m$, are independent Poisson processes with intensity $\lambda_j$ taking values in $\{mh_j,\ m\in\mathbb{N}\}$. The path is still deterministic with drift $b'$ and has jumps of size in $\{|h_1|,\dots, |h_m|\}$ occurring at exponentially distributed random times. When we let $m$ tend to $\infty$ in a suitable sense, or more generally when the L\'evy measure $\nu$ is of finite variation, that is, condition \eqref{eq:finite_variation} holds, we can write
$$
X(t)=b't+\sum_{0\leq s\leq t}\Delta X(s),
$$
where $\Delta X(s) = X(s) - X(s-)$ is the jump at time $s$. Instead of dealing with jumps directly, it is more convenient to count the number of jumps that belong to a set $A$ up to time $t$. To this end, for a Borel set $A\subseteq\mathbb{R}^n\setminus\{0\}$ and $t\geq 0$, we define the random Poisson measure with intensity $\nu$
$$
N(t,A)=\#\{0\leq s\leq t\ |\ \Delta X(s)\in A\},
$$
which allows us to write
$$
\sum_{0\leq s\leq t}\Delta X(s)=\int_{\mathbb{R}^n\setminus\{0\}} xN(t,dx).
$$
However, in the most general case, the L\'evy measure $\mu$ may not satisfy the finite variation condition \eqref{eq:finite_variation} and to deal with the accumulation of small jumps, we make use of the compensated Poisson measure:
$$
\widetilde{N}(dt, dx) = N(dt, dx) - dt\ \nu(dx).
$$
%inserted space in second term
Finally, in case of a general L\'evy measure $\nu$ and of a diffusion matrix $A$, one has the \emph{L\'evy-It\^o decomposition} \cite[Theorem 2.4.16]{Applebaum}:
\begin{equation}
\label{eq:Generator}
X(t)= D W(t) + bt+\int_{0<|x|<1}x\widetilde{N}(t, dx) + \int_{|x|\geq 1}xN(t,dx),
\end{equation}
where $D$ is a $n\times n$-dimensional matrix such that $DD^T=A$, and $\{W(t)\}_{t\geq 0}$ is a $n$-dimensional Brownian motion.

\subsection{Connections to integro-differential operators}
\label{sec:Connection_to_pde}
At this point we want to explore the connection between stochastic processes and integro-differential operators. Using the fact that any L\'evy process is a Markov process, by defining
$$
T_tf(x) := \EE\left[f(x+X(t))\right],\quad\forall\,x\in\RR^n,
$$
we obtain that $\{T_t\}_{t\geq 0}$ defines a one-parameter semigroup of linear operators on the Banach space of bounded continuous functions, $C(\mathbb{R}^n)$. One can think of the semigroup $\{T_t\}_{t\geq 0}$ as a tool to give a deterministic, macroscopic description of the L\'evy process as an average of microscopic random dynamics. The infinitesimal generator corresponding to the semigroup semigroup $\{T_t\}_{t\geq 0}$ is defined formally by
$$
Lf(x) = \lim_{t\downarrow 0}\frac{T_tf(x) - f(x)}{t},
$$
and takes the form
\begin{align*}
Lf(x)&=\frac{1}{2} \mbox{tr}(AD^2f) + b\cdot\nabla f(x)
+\int_{\mathbb{R}^n\setminus\{0\}}\left[f(x+y)-f(x)-y\cdot\nabla f(x)\chi_{|y|<1}(y)\right]\ \nu(dy).
\end{align*}
Under suitable regularity assumptions that allow us to apply It\^o's rule \cite[Theorem 4.4.7]{Applebaum} to solutions to the parabolic differential equation $u_t = Lu$ on $(0,\infty)\times\RR^n$, with initial condition $u(0, \cdot) = f$ on $\RR^n$, we obtain that $u(t, x) = T_tf(x)$, for all $t \geq 0$ and $x\in \RR^n$, and so $T_t = e^{tL}$.

We can also establish a connection between the infinitesimal generator $L$ of the process $\{X(t)\}_{t\geq 0}$ and the characteristic exponent $\psi(\xi)$ appearing in the L\'evy-Khintchine formula \eqref{eq:LK}. Viewed as a pseudo-differential operator \cite{Boyarchenko_Levendorskii_2002b, Taylor_vol2}, the symbol of the infinitesimal generator $L$ is the characteristic exponent \eqref{eq:exponent} appearing in identity \eqref{eq:LK}. In our survey, we will be concerned with generalizations of symbols that contain only a drift and a nonlocal term (the second order diffusion term is removed). This gives rise to mathematical challenges in the study of the regularity of solutions, when the drift term dominates the nonlocal component -- the so-called \emph{supercritical regime}. This property is often encountered in financial models for stock prices, such as Variance Gamma and Regular L\'evy Processes of Exponential Type described in greater detail in \S\ref{examples}.

\subsection{Stochastic integro-differential equations}
\label{sec:SDE}
More generally than the infinitesimal generators of L\'evy processes, in this survey we are specifically concerned with \emph{nonlocal} operators that are infinitesimal generators of strong Markov processes, which can be written as solutions to stochastic integro-differential equations of the form:
\begin{equation}
\label{eq:SDE}
dX(t)=b(X(t-))dt+\int_{\mathbb{R}^n\setminus\{0\}}F(X(t-),y)\tilde{N}(dt,dy),\quad t>0.
\end{equation}
Here $\tilde{N}(dt,dy)$ is a compensated Poisson random measure with intensity measure $d\nu$, as defined in \S\ref{sec:Rep_Levy_proc},
and $b$ and $F$ satisfy suitable conditions, which we describe in detail in \S\ref{Statements}. Our conditions ensure, by \cite[Theorem 6.2.9]{Applebaum}, that for any initial condition $X^x(0)=x\in\mathbb{R}^n$, there exists a unique strong solution $\{X^x(t)\}_{t\geq 0}$ to equation \eqref{eq:SDE} with \emph{c\`adl\`ag} paths a.s. The process $\{X^x(t)\}_{t\geq 0}$ satisfies the strong Markov property, and therefore it is uniquely determined by its infinitesimal generator
\begin{align}
\label{main}
Lu(x)=b\cdot\nabla u(x)+\int_{\mathbb{R}^n\setminus\{0\}}\left(u(x+F(x,y))-u(x)-F(x,y)\cdot\nabla u(x)\right)\ \nu(dy)
\end{align}
for all $u\in C^2(\mathbb{R}^n)$ (this denotes all functions with bounded and continuous derivatives up to and including order 2 in $\mathbb{R}^n)$. The term \emph{nonlocal} refers to the fact that the value of $Lu(x)$ depends on the whole solution $u$ and not only on
its behavior nearby the point $x$.
% Maybe this example can be removed because we don't use in in our results
A typical example of a nonlocal integro-differential operator is the fractional Laplacian $(-\Delta)^{s}$, with $s\in (0,1)$, which is defined on the Fourier transform side by the formula
$$
\widehat{(-\Delta)^{s}u}(\xi)=|\xi|^{2s}\hat{u}(\xi),
$$
or, equivalently, by the pointwise representation
$$
{(-\Delta)^{s}u}(x)=\gamma(n,s)\ p.v.\int_{\mathbb{R}^n}\frac{2u(x)-u(x+y)-u(x-y)}{|x|^{n+2s}}\ dy,
$$
$\gamma$ being a normalization constant depending only on $n$ and $s$. The fractional Laplacian $(-\Delta)^s$
%removed comma
is the infinitesimal generator of the symmetric $2s$\emph{-stable L\'evy process} with characteristic exponent in the L\'evy-Khintchine formula given by $\psi(\xi) = |\xi|^{2s}$.

\subsection{Obstacle problems}
\label{sec:Obs_prob}
\begin{figure}
\begin{center}
\begin{tikzpicture}

\draw [black, thick, fill=lightgray]   plot[smooth, tension=.7] coordinates
{(-5.5,0) (-5,2) (-3,3) (-1.5,2.5) (2,2.5) (3,2) (3,0.5) (0.5,-1) (-2,-0.5) (-5,-1) (-5.5,0)};

\draw[draw=black,thick,<-] (3,0.4) .. controls (4,0) .. (4.5,-0.5);
\node[right] at (4.5,-0.5) {exercise (free) boundary};

%exercise region
\node[below] at (-2,2) {$\{v=\varphi\}$};
\node[below] at (-2,1) {$-\partial_t v-Lv+rv \geq 0$};

%continuation region
\node[below] at (6,2) {$\{v>\varphi\}$};
\node[below] at (6,1) {$-\partial_t v-Lv+rv = 0$};

\end{tikzpicture}
\caption{A schematic description of the complementarity conditions for the evolution obstacle problem at a time slice $t$. The exercise region $\{v=\varphi\}$ is represented by the gray area, and the remaining region is the continuation region $\{v>\varphi\}$.}
\label{fig:Obs_prob_evol}
\end{center}
\end{figure}
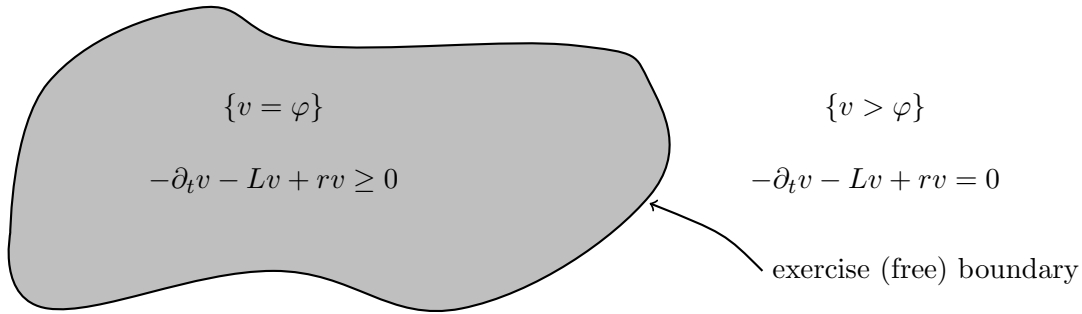

In recent years there has been a resurgence of interest in the study of nonlocal operators, motivated by applications. In fact, such operators and the associated integro-differential equations naturally arise in a variety of contexts, ranging from temperature control to linear elasticity, from fluid dynamics to financial mathematics. To describe the latter application in more detail, we assume that
\begin{equation}\label{asset}
S(t)=e^{X(t)}
\end{equation}
models an asset price process, where $\{X(t)\}_{t\geq 0}$ is a solution to the stochastic equation \eqref{eq:SDE}. We let ${\varphi:\mathbb{R}^n\to \mathbb{R}}$ be the payoff function of an American option (i.e., a profit of $\varphi(s)$ is generated when exercising the option at time $t$ and the stock level is $s=S(t)$). Without loss of generality, we can assume that the payoff can be written as a function of $\{X(t)\}_{t\geq 0}$. We recall that, unlike the European option, in the American option framework the holder has the right to exercise at any date prior to maturity, and not only at the expiry date. Hence, the value of the American option with expiry date $T$ can be written as
$$
v(t, x)=\sup\ \mathbb{E}[e^{-rt}\varphi (X(\theta)) | X(t) = x], \quad\hbox{ for all } (t, x) \in (0,T)\times\RR^n,
$$
where the supremum is taken over all stopping times $\theta$ bounded by $T-t$, and we assume that the expectation is taken under a risk-neutral probability measure and $r$ is the risk-free interest rate. Letting $\tau$ be the first time that the stochastic process $\{X(t)\}_{t\geq 0}$ enters the \emph{exercise region } ${\{v=\varphi\}}$, and assuming that the value function $u(t,x)$ is regular enough, probabilistic arguments ensure that the stopped process $\{e^{-rt\wedge\tau}v(t\wedge\tau,X(t\wedge\tau))\}_{t\geq 0}$ is a martingale, which is equivalent to the equality
\begin{equation}
\label{eq:Martingale}
\partial_t v + Lv - rv = 0, \quad\hbox{ for all } (t, x) \in \{v>\varphi\}.
\end{equation}
In general, however the discounted option price process $\{e^{-rt}v(t,X(t))\}_{t\geq 0}$ is a supermartingale, which translates into the inequality
\begin{equation}
\label{eq:Supermartingale}
\partial_t v + Lv - rv \leq 0, \quad\hbox{ for all } (t, x) \in (0,T)\times\RR^n.
\end{equation}
Combining equations \eqref{eq:Martingale} and \eqref{eq:Supermartingale} together with the property that $v\geq \varphi$ gives us that the value function $v$ is  a solution to the \emph{evolution obstacle problem}:
\begin{equation}
\label{obstacle}
\min\{-\partial_t v-Lv+rv, v-\varphi\}=0, \quad\hbox{ for all } (t, x) \in (0,T)\times\RR^n,
\end{equation}
where $L$ is the infinitesimal generator of $\{X(t)\}_{t\geq 0}$. The strong Markov property of $\{X(t)\}_{t\geq 0}$ implies that the exercise decision at any time $t$ depends only on $t$ and $X(t)$. Therefore, for each $t$ there exist an  \emph{exercise region } ${\{v=\varphi\}}$, in which one should exercise the option,  and a { \emph{continuation region }} ${\{v>\varphi\}}$, in which one instead should wait. The {\emph{exercise boundary}} is the interface separating the two. See Figure \ref{fig:Obs_prob_evol} for a schematic representation. We briefly mention here that in the case of perpetual American option, when the option has a infinite expiration time, the value function depends only on the current value of the process $\{X(t)\}_{t\geq 0}$ and is a solution to a stationary obstacle problem. We refer to \S\ref{Statements} for further details.

\subsection{Review of literature and outline of the survey}
\label{sec:Literature_overview}
If the underlying stochastic process is {Brownian motion}, then the infinitesimal generator of the underlying process is ${L=\Delta}$ and $u$ will satisfy the classical obstacle problem, which is by now very well understood \cite{Caffarelli_1977, Caffarelli_1980, Caffarelli_jfa_1998, Laurence_Salsa_2009}.  However, Brownian motion falls short in some respects:
\begin{itemize}
\item[1.] Stock prices do not move continuously, which prompts us to consider models that allow jumps in small time intervals;
\item[2.] Empirical studies of stock price returns indicate distributions with heavy tails, which are not compatible with a Gaussian model.
\end{itemize}
For these reasons, it becomes necessary to  study  jump-diffusion processes, whose infinitesimal generator is an integro-differential operator of the form \eqref{eq:Generator}. Such type of operators were introduced in finance by the Nobel Prize winner R. Merton \cite{Merton1976}. The novel element, which reflects the presence of jumps, is the integral term.  Its presence leads to new theoretical and numerical issues. Since no closed form solutions are known in general for the American option,  it becomes important to determine the regularity of the exercise boundary, which in turn is closely related to the behavior of the value function.

In the framework of jump-diffusion models with a non-degenerate diffusion matrix, regularity of the value function and efficient numerical schemes were studied in \cite{Levendorski_2008, Bayraktar_2009, Bayraktar_2011, Bayraktar_Xing_2012}, and regularity of the free boundary was explored in \cite{Bayraktar_Xing_2009sjma}. Using methods from the theory of pseudo-differential operators and the Wiener-Hopf factorization, qualitative studies of American option prices and of the exercise region under pure-jump models were performed in articles such as \cite{Mordecki_2002, Levendorski_2004a, Levendorski_2004b, Boyarchenko_Levendorskii_2011, Boyarchenko_Levendorskii_2002b}.

Our work continues the study of the regularity of solutions to obstacle problems for nonlocal operators with (possibly supercritical) drift. The purpose of this note is to give an overview of the regularity results obtained in \cite{Danielli_Petrosyan_Pop}. In \S2, we describe two examples of stochastic processes of interest in mathematical finance to which our results apply. In \S3, we state the problem precisely, and provide the statements of our main results.

\section{Motivating Examples}
\label{examples}

In this section we assume $n=1$  and that the asset price process can be written as in \eqref{asset}. Moreover, $r$ denotes the risk-free interest rate.  It is crucially important to ensure that the discounted asset price process $\{e^{-rt} S(t)\}_{t\geq 0}$ is a martingale in order to obtain an arbitrage-free market. Assume that $\{X(t)\}_{t\geq 0}$ is a one-dimensional L\'evy process that satisfies the stochastic equation:
\begin{equation}
\label{eq:Process_const_drift}
dX(t) = b\, dt + \int_{\mathbb{R}^n} y \, \widetilde N(dt, dy),\quad\forall\, t>0,
\end{equation}
where $b$ is a real constant and $\widetilde N(dt, dy)$ is a compensated Poisson random measure with L\'evy measure $\nu(dy)$. Using \cite[Theorem~5.2.4 and Corollary~5.2.2]{Applebaum},  a sufficient
condition that guarantees that the discounted asset price process $\{e^{-rt+X(t)}\}_{t\geq 0}$ is a martingale is:
\begin{equation}
\label{eq:Martingale_cond}
\begin{aligned}
\int_{|x|\geq 1} e^x\, \nu(dx) &<\infty
\quad\hbox{ and }\quad
-r+\psi(-i) &= 0,
\end{aligned}
\end{equation}
where $\psi(\xi)$ denotes the characteristic exponent of the L\'evy process $\{X(t)\}_{t\geq 0}$, that is,
\begin{equation}
\label{eq:Characteristic_exponent}
\psi(\xi) = ib\xi+\int_{\mathbb{R}\setminus\{0\}} (e^{ix\xi}-1-ix\xi)\,\nu(dx).
\end{equation}
Examples in mathematical finance to which our results
apply include the \emph{Variance Gamma Process} \cite{Madan_Seneta_1990} and
\emph{Regular L\'evy Processes of Exponential type} (RLPE) \cite{Boyarchenko_Levendorskii_2002b}.

When the jump-part of the nonlocal operator $L$ corresponding to the integral term in the characteristic exponent \eqref{eq:Characteristic_exponent} has sublinear growth as $|\xi|\rightarrow\infty$, we say that the drift term $b\cdot\nabla$ corresponding to $ib\cdot\xi$ in the characteristic exponent \eqref{eq:Characteristic_exponent} is \emph{supercritical}. An example of a nonlocal operator with supercritical drift is the Variance Gamma Process and a subcollection of Regular L\'evy Processes of Exponential type described below.

\subsection{Variance Gamma Process}
\label{sec:VG}
Following \cite[Identity (6)]{Carr_Geman_Madan_Yor_2002}, the Variance Gamma Process $\{ X(t)\}_{t\geq 0}$ with parameters $\nu,\sigma,$ and $\theta$ has L\'evy measure given by
\begin{align*}
\nu(dx) = \frac{1}{\nu|x|}\left(e^{-\frac{|x|}{\eta_p}}\mathbf{1}_{\{x>0\}} + e^{-\frac{|x|}{\eta_n}}\mathbf{1}_{\{x<0\}}\right)\, dx,
\end{align*}
where $\eta_p>\eta_n$ are the roots of the equation $x^2-\theta\nu x-\sigma^2\nu/2=0$, and $\nu,\sigma, \theta$ are positive constants. From \cite[Identity (4)]{Carr_Geman_Madan_Yor_2002}, we have that the characteristic exponent of the Variance Gamma Process with constant drift $b\in\mathbb{R}$, $\{X(t)+bt\}_{t\geq 0}$, has the expression:
$$
\psi_{\hbox{\tiny{VG}}}(\xi) = \frac{1}{\nu}\ln \left(1-i\theta\nu\xi + \frac{1}{2}\sigma^2\nu\xi^2\right) + ib\xi,
\quad\forall\, \xi\in\mathbb{C},
$$
and so the infinitesimal generator of $\{X(t)+bt\}_{t\geq 0}$ is given by
$$
L = \frac{1}{\nu}\ln(1-\theta\nu\nabla - \frac{1}{2}\sigma^2\nu\Delta)+b\cdot\nabla,
$$
which is a sum of a pseudo-differential operator of order less
than any $s>0$ and one of order $1$. When $\eta_p<1$ and $r=\psi_{VG}(-i)$, condition \eqref{eq:Martingale_cond} is satisfied and the discounted asset price process $\{e^{-rt+X(t)}\}_{t\geq 0}$ is a martingale. Thus, applying the results in \S\ref{Statements} to the Variance Gamma Process $\{X(t)\}_{t\geq 0}$ with constant drift $b$, we obtain that the prices of perpetual and finite expiry American options with bounded and Lipschitz payoffs are Lipschitz functions in the spatial variable. Given that the nonlocal component of the infinitesimal generator $L$ has order less than any $s>0$, this may be the optimal regularity of solutions that we can expect.

\subsection{Regular L\'evy Processes of Exponential type}
\label{sec:RLPE}
Following \cite[Chapter~3]{Boyarchenko_Levendorskii_2002b}, for parameters $\lambda_{-} < 0 < \lambda_{+}$, a L\'evy process is said to be of exponential type $[\lambda_{-}, \lambda_{+}]$ if it has a L\'evy measure $\nu(dx)$ such that
$$
\int_{-\infty}^{-1} e^{-\lambda_{+}x} \nu(dx) + \int_1^{\infty} e^{-\lambda_{-}x}\nu(dx) <\infty.
$$
Regular L\'evy Processes of Exponential type $[\lambda_{-}, \lambda_{+}]$ and order $\nu$ are non-Gaussian L\'evy
processes of exponential type $[\lambda_{-}, \lambda_{+}]$ such that, in a neighborhood of zero, the L\'evy measure can be represented as $\nu(dx) = f(x)\,dx$, where the density $f(x)$ satisfies the property that
$$
|f(x) - c|x|^{-\nu-1}| \leq C|x|^{-\nu'-1},\quad\forall\, |x|\leq 1,
$$
for constants $\nu'<\nu$, $c>0$, and $C>0$. Our results apply to RLPE type $[\lambda_{-}, \lambda_{+}]$, when we choose the parameters $\lambda_-\leq-1$ and $\lambda_+ \geq 1$. The class of RLPE include the CGMY/KoBoL processes introduced in \cite{Carr_Geman_Madan_Yor_2002}. Following \cite[Equation (7)]{Carr_Geman_Madan_Yor_2002}, CGMY/KoBoL processes are characterized by a L\'evy measure of the form
$$
\nu(dx) = \frac{C}{|x|^{1+Y}}\left(e^{-G|x|}\mathbf{1}_{\{x<0\}} + e^{-M|x|}\mathbf{1}_{\{x>0\}}\right)\, dx,
$$
where the parameters $C>0$, $G, M \geq 0$, and $Y<2$. Our results
apply to CGMY/KoBoL processes, when we choose the parameters $G,M>1$ and
$Y<2$, or $G,M\geq 1$ and $0<Y<2$.

\section{Statements of the main results}
\label{Statements}
In this section we provide the statements of our main results. Complete
proofs can be found in \cite{Danielli_Petrosyan_Pop}, where these
results have originally appeared.

We begin by listing the required assumptions on the measure $\nu(dx)$ and the coefficients $b(x)$ and $F(x,y)$ appearing in the operator \eqref{main}:
\begin{enumerate}
\item[1.] There is a positive constant $K$ such that for all $x_1, x_2\in\mathbb{R}^n$, we have
\begin{equation*}
\begin{aligned}
&\int_{\mathbb{R}^n\setminus\{O\}} |F(x_1, y)-F(x_2, y)|^2 \, d\nu(y) \leq K |x_1-x_2|^2,\\
&\sup_{z\in B_{|y|}}|F(x,z)| \leq \rho(y),\quad\forall\, x,y\in\mathbb{R}^n,\\
&\int_{\mathbb{R}^n\setminus\{O\}} \left(|y|\vee\rho(y)\right)^2 \,\nu(dy) \leq K,
\end{aligned}
\end{equation*}
where $\rho:\mathbb{R}^n\rightarrow[0,\infty)$ is a measurable function.
\item[2.]
The coefficient $b:\mathbb{R}^n\rightarrow\mathbb{R}^n$ is bounded and Lipschitz continuous, i.e., $b\in C^{0,1}(\mathbb{R}^n)$.
\item[3.] For the stationary problem, we assume that $F(x, y) = F(y)$
(independent of $x$).
\end{enumerate}

\subsection{Stationary obstacle problem}
\label{sec:Stationary}
We consider the obstacle problem
\begin{equation}\label{stationary}
\min\{-Lv + c v - f, v - \varphi\} = 0\quad\hbox{ on }\mathbb{R}^n,
\end{equation}
where $L$ is the infinitesimal generator of the unique strong solution $\{X^x(t)\}_{t\geq 0}$ to the stochastic equation \eqref{eq:SDE}, with initial condition $X^x(0)=x$. We explicitly remark here that, in the applications in \S\ref{examples}, one chooses $c\equiv r$, the risk-free interest rate.
Solutions to the obstacle problem \eqref{stationary} are constructed using the \emph{{stochastic representation formula}} of the value function:
\begin{equation*}
{v(x) := \sup\{v(x;\tau):\,\tau\in\mathcal{T}\}}.
\end{equation*}
where $\cT$ is the set of stopping times and
\begin{equation*}
{v(x;\tau) :=\mathbb{E} \left[e^{-\int_0^{\tau} c(X^x(s))\, ds} \varphi(X^x(\tau))
+ \int_0^{\tau} e^{-\int_0^ t c(X^x(s))\, ds} f(X^x(t))\, dt\right]},\quad\forall\,\tau\in\cT.
\end{equation*}
In order to state our results, we need to introduce the relevant function spaces. We denote by $C(\mathbb{R}^n)$ the space of bounded continuous functions $u:\mathbb{R}^n\rightarrow\mathbb{R}$ such that
\begin{equation*}
\|u\|_{C(\mathbb{R}^n)} :=\sup_{x\in \mathbb{R}^n} |u(x)| < \infty.
\end{equation*}
For all $\alpha\in (0,1]$, a function $u: \mathbb{R}^n\rightarrow\mathbb{R}$ belongs to the H\"older space of functions $C^{0,\alpha}(\mathbb{R}^n)$ if
\begin{align*}
\|u\|_{C^{0,\alpha}(\mathbb{R}^n)} := \|u\|_{C(\mathbb{R}^n)} + [u]_{C^{0,\alpha}(\mathbb{R}^n)}<\infty,
\end{align*}
where, as usual, we define
\begin{equation*}
[u]_{C^{0,\alpha}(\mathbb{R}^n)}:= \sup_{x_1,x_2\in\mathbb{R}^n, x_1\neq x_2} \frac{|u(x_1)-u(x_2)|}{|x_1-x_2|^{\alpha}}.
\end{equation*}
When $\alpha\in (0,1)$, we denote for brevity $C^{\alpha}(\mathbb{R}^n) := C^{0,\alpha}(\mathbb{R}^n)$. Our first result concerns the regularity of the value function.

\begin{thrm}\label{VF} Let $c,\varphi,f:\mathbb{R}^n\rightarrow \mathbb{R}$ be bounded Lipschitz continuous functions, and assume that there is a constant $c_0>0$ such that ${c(x) \geq c_0>0,\ \forall\, x\in\mathbb{R}^n}$.
{Then the following hold:}
\begin{enumerate}
\item[(i)] \emph{{(H\"older continuity)}}
{There is a constant $\alpha = \alpha([b]_{C^{0,1}(\mathbb{R}^n)}, c_0) \in (0,1)$, such that the value function} {$v\in C^{\alpha}(\mathbb{R}^n)$}.
\item[(ii)] \emph{{(Lipschitz continuity)}}
{If in addition we have that}
\begin{equation}
\label{eq:Zeroth_oreder_term_cond}
{c_0\geq [b]_{C^{0,1}(\mathbb{R}^n)}},
\end{equation}
{then the value function} ${v\in C^{0,1}(\mathbb{R}^n)}$.
\end{enumerate}
\end{thrm}

The proof of Theorem \ref{VF} hinges on the
stochastic representation of solutions and on the continuity of
the strong solutions to the SDE with respect to the initial
conditions.
To proceed, we introduce the notion of viscosity solution, which gives an intrinsic definition of a solution which is local in nature, but does not assume  a priori any regularity, except for continuity.

\begin{dfn}
Let $v\in C(\mathbb{R}^n)$. We say that $v$ is a \emph{viscosity subsolution (supersolution)} to the stationary obstacle problem
 if, for all $u\in C^2(\mathbb{R}^n)$ such that $v-u$ has a global max (min) at $x_0\in\mathbb{R}^n$ and
$u(x_0)=v(x_0)$, then
\begin{equation}
{\min\{-Lu(x_0) + c(x_0) u(x_0) - f(x_0), u(x_0) - \varphi(x_0)\} \leq (\geq)\, 0}.
\end{equation}
{We say that $v$ is a \emph{viscosity solution}  if it is both a sub- and supersolution.}
\end{dfn}
Next, we show that the value function is the unique solution to \eqref{stationary}.
\begin{thrm}[Existence]
{Assume in addition}
\begin{equation*}
{\int_{\mathbb{R}^n\setminus\{O\}} |F(y)|^{2\alpha}\, \nu(dy) <\infty}
\end{equation*}
{where $\alpha\in (0,1)$ is the constant appearing in Theorem 1.
Then the value function $v$ is a viscosity solution to the stationary obstacle problem.}
\end{thrm}

\begin{thrm}[Uniqueness]
{Suppose that  $c,f,\varphi\in C(\mathbb{R}^n)$ and $c$ is a positive function}.
If the stationary obstacle problem  has a viscosity solution, then it is unique.
\end{thrm}

We remark that a sufficient condition on the L\'evy measure to ensure that perpetual American put option prices are Lipschitz continuous, but not continuously differentiable, is provided in \cite[Theorem~5.4, p.~133]{Boyarchenko_Levendorskii_2002b}. However, the condition is in terms of the Wiener-Hopf factorization for the characteristic exponent of the L\'evy process, and it is difficult to find a concrete example for which it holds.
Since in our case the order of the nonlocal operator is strictly
less than the order of the drift component, and there is no second-order term, the issue of regularity of  solutions is quite delicate.

The proof of the existence result hinges in a crucial way on a \emph{Dynamic Programming Principle}. In order to state it precisely, we need the following definition.
\begin{dfn}
For all $r>0$ and $x\in\mathbb{R}^n$, we let
\begin{equation*}
{\tau_r := \inf\{t \geq 0: X^x(t)\notin B_r(x)\}},
\end{equation*}
where $B_r(x)$ denoted the open Euclidean ball of radius $r>0$ centered at $x\in\RR^n$.
\end{dfn}

\begin{thrm}[Dynamic Programming Principle]
\label{lem:DPP_stat}
{The value function $v(x)$ satisfies:}
\begin{equation*}
{v(x) = \sup\{v(x;r,\tau):\,\tau \leq \tau_r\}},\quad{\forall\, r>0},
\end{equation*}
{where we define}
\begin{equation*}
\begin{aligned}
{v(x;r,\tau)} &{:=\mathbb{E}\left[e^{-\int_0^{\tau}c(X^x(s))\,ds} \left(\varphi(X^x(\tau))\mathbf{1}_{\{\tau<\tau_r\}}+ v(X^x(\tau))\mathbf{1}_{\{\tau =\tau_r\}}\right)\right]}\\
&\quad  {+\mathbb{E}\left[\int_0^{\tau\wedge\tau_r} e^{-\int_0^t c(X^x(s))\, ds} f(X^x(t))\, dt\right]}.
\end{aligned}
\end{equation*}
\end{thrm}

Uniqueness is proved instead  with the aid of the following theorem.

\begin{thrm}[Comparison principle]
{Suppose that the assumptions of the uniqueness theorem hold.
 If $u$ and $v$ are a viscosity subsolution and supersolution to the stationary obstacle problem, respectively, then $u\leq v$.}
\end{thrm}
In financial terms, comparison principles simply translate into arbitrage inequalities:
if the terminal payoff of an American option dominates the terminal
payoff of another one, then their values should verify the same inequality.

\subsection{Evolution Obstacle Problem}
The evolution obstacle problem is given by
\begin{equation}
\label{evol}
\begin{cases}
{\min\{-\partial_t v-Lv + c v - f, v - \varphi\} = 0}& \quad\hbox{ on } [0,T)\times\mathbb{R}^n,\\
{v(T,\cdot) = g}&\quad\hbox{ on } \mathbb{R}^n,
\end{cases}
\end{equation}
with the compatibility condition
\begin{equation}
\label{eq:Compatibility_evol}
{g\geq \varphi(T,\cdot)\quad\hbox{ on }\quad \mathbb{R}^n}.
\end{equation}
The treatment of this problem is very similar to the stationary case.  For the sake of brevity, we confine ourselves to mentioning here that the main new difficulty is to establish regularity in the time variable.
This is done with the aid of the following result concerning the continuity properties of $\{X(t)\}_{t\geq 0}$, which in turn is a consequence of Doob's Martingale Inequality.

%We could remove this lemma
\begin{lemma}
{There is a positive constant $C=C(\|b\|_{C^{0,1}(\mathbb{R}^n)}, K)$ such that}
\begin{align*}
{\mathbb{E}\left[\max_{s\in [0,t]} \left|X^{x_1}(s)-X^{x_2}(s)\right|^2\right] }&{\leq C|x_1-x_2|^2 e^{Ct},
\quad\forall\, x_1,x_2\in\mathbb{R}^n,\, t\geq 0,}\\
{\mathbb{E}\left[\max_{r\in [s,t]} \left|X^x(r)-X^x(s)\right|^2\right]} &{\leq C|t-s|\vee|t-s|^2,\quad\forall\, x\in\mathbb{R}^n,\, 0\leq s<t.}
\end{align*}
\end{lemma}
The use of this lemma also allows  to relax the assumptions on the coefficients in that we no longer require condition \eqref{eq:Zeroth_oreder_term_cond} to hold and we can allow the jump size $F(x, y)$ to be a function of the current state $x$ of the process.

The relevant function spaces, in the evolution case, are as follows. For all $T>0$, we denote by $C^{\frac{1}{2}}_tC^{0,1}_x([0,T]\times\mathbb{R}^n)$ the space of functions $u:[0,T]\times\mathbb{R}^n\rightarrow\mathbb{R}$ such that
\begin{equation*}
\|u\|_{C^{\frac{1}{2}}_tC^{0,1}_x([0,T]\times\mathbb{R}^n)}:= \|u\|_{C([0,T]\times\mathbb{R}^n)}
+
\sup_{\stackrel{t_1,t_2\in [0,T], t_1\neq t_2}{x_1,x_2\in\mathbb{R}^n, x_1\neq x_2}}
\frac{|u(t_1,x_1)-u(t_2,x_2)|}{|t_1-t_2|^{\frac{1}{2}}+|x_1-x_2|} <\infty,
\end{equation*}
and we let $C^1_tC^2_x([0,T]\times\mathbb{R}^n)$ denote the space of functions $u:[0,T]\times\mathbb{R}^n\rightarrow\mathbb{R}$ such that the first order derivative in the time variable and the second order derivatives in the spatial variables are continuous and bounded. Let $\mathcal{T}_t$ denote the set of stopping times $\tau\in\mathcal{T}$ bounded by $t$, for all $t\geq 0$. Solutions to problem \eqref{evol} are constructed using the stochastic representation formula,
\begin{equation}
\label{eq:Value_function_evol}
v(t,x) := \sup \{ v(t,x;\tau):\, \tau\in\mathcal{T}_{T-t}\},
\end{equation}
where we define
\begin{equation}
\label{eq:Value_function_evol_aux}
\begin{aligned}
v(t,x;\tau) &:= \mathbb{E}\left[e^{-\int_0^{\tau} c(t+s, X^x(s))\, ds}\varphi(t+\tau, X^x(\tau)) \mathbf{1}_{\{\tau<T-t\}}\right]\\
&\quad+ \mathbb{E} \left[e^{-\int_0^{\tau} c(t+s, X^x(s))\, ds}g(X^x(T-t)) \mathbf{1}_{\{\tau = T-t\}}\right]\\
&\quad+ \mathbb{E} \left[\int_0^{\tau} e^{-\int_0^s c(t+r,X^x(r))\, dr} f(t+s,X^x(s))\, ds\right],
\end{aligned}
\end{equation}
for all $(t,x)\in[0,T]\times\in\mathbb{R}^n$.

\begin{prp}[Regularity]
\label{prop:Regularity_evol}
Suppose that $c,\varphi,f$ belong to $C^{0,1}([0,T]\times\mathbb{R}^n)$, the final condition $g$ is in $C^{0,1}(\mathbb{R}^n)$, and the compatibility condition \eqref{eq:Compatibility_evol} holds. Then the value function $v$ defined in \eqref{evol} belongs to $C^{\frac{1}{2}}_tC^{0,1}_x([0,T]\times\mathbb{R}^n)$.
\end{prp}

We next define a notion of viscosity solution for the evolution obstacle problem \eqref{evol} extending that of its stationary analogue for equation \eqref{stationary} similarly to the ideas described in \cite[\S8]{Crandall_Ishii_Lions_1992}:

\begin{dfn}[Viscosity solutions]
\label{defn:Viscosity_sol_evol}
Let $v\in C(\RR^n)$. We say that $v$ is a viscosity subsolution (supersolution) to the evolution obstacle problem
\eqref{evol} if
\begin{equation}
\label{eq:Final_cond_evol}
v(T,\cdot) \leq (\geq) g,
\end{equation}
and, for all $u\in C^1_tC^2_x([0,T]\times\RR^n)$ such that
$v-u$ has a global max (min) at $(t_0,x_0)\in[0,T)\times\RR^n$ and $u(t_0,x_0)=v(t_0,x_0)$, we have that
\begin{equation}
\label{eq:Sub_super_sol_evol}
\min\{-\partial_tu(t_0,x_0)-Lu(t_0,x_0) + c(t_0,x_0) u(t_0,x_0) - f(t_0,x_0), u(t_0,x_0) - \varphi(t_0,x_0)\} \leq (\geq)\, 0.
\end{equation}
We say that $v$ is a viscosity solution to equation \eqref{evol} if it is both a sub- and supersolution.
\end{dfn}

\begin{thrm}[Existence]
\label{thm:Existence_evol}
Suppose that the hypotheses of Proposition~\ref{prop:Regularity_evol} hold. Then the value function $v$ defined in \eqref{eq:Value_function_evol} is a viscosity solution to the evolution obstacle problem \eqref{evol}.
\end{thrm}

We conclude with

\begin{thrm}[Uniqueness]
\label{thm:Uniqueness_evol}
Suppose that $g$ belongs to $C(\RR^n)$, $c,f,\varphi$ are in $C([0,T]\times\RR^n)$, the compatibility condition \eqref{eq:Compatibility_evol} holds, and
\begin{equation}
\label{eq:Interior_point}
\lim_{y\to O} F(x,y) = 0,\quad\forall\, x\in\RR^n.
\end{equation}
If the obstacle problem \eqref{evol} has a solution, then it is unique.
\end{thrm}

%%%%%%%%%%%%%%%%%%%%%%%%%%%%%%%%%%%%%%%%%%%%%%%%%%%%%%%%%%%%%%%%%%%%%%%%%%%%%%%
%
%                                bibliography
%
%%%%%%%%%%%%%%%%%%%%%%%%%%%%%%%%%%%%%%%%%%%%%%%%%%%%%%%%%%%%%%%%%%%%%%%%%%%%%%%

\bibliography{mfpde}
\bibliographystyle{amsplain}

\end{document}

%% file: latexformat.tex
\newtheorem*{thm*}{Theorem}

\newtheorem*{lem*}{Lemma}

\theoremstyle{definition}

\renewcommand{\thecase}{}
\newtheorem*{case*}{Case}

\newtheorem*{defn*}{Definition}

\newtheorem*{exmp*}{Example}

\newtheorem*{rmk*}{Remark}
\renewcommand{\thestep}{}

\theoremstyle{remark}

\makeatletter
\def\alphenumi{
  \def\theenumi{\alph{enumi}}
  \def\p@enumi{\theenumi}
  \def\labelenumi{(\@alph\c@enumi)}}
\makeatother

\makeatletter
\def\thecase{\@arabic\c@case}
\makeatother

\makeatletter
\def\thestep{\@arabic\c@step}
\makeatother

%% file: latexnewmacros.tex
\newcount\hh
\newcount\mm
\mm=\time
\hh=\time
\divide\hh by 60
\divide\mm by 60
\multiply\mm by 60
\mm=-\mm
\advance\mm by \time
\def\hhmm{\number\hh:\ifnum\mm<10{}0\fi\number\mm}

\let\oldmarginpar\marginpar
\renewcommand\marginpar[1]{\-\oldmarginpar[\raggedleft\footnotesize #1]%
{\raggedright\footnotesize #1}}

\newcommand\EE{\mathbb{E}}

\newcommand\PP{\mathbb{P}}

\newcommand\RR{\mathbb{R}}

\newcommand\cT{{\mathcal{T}}}

\newcommand\eps{\varepsilon}